
\magnification=1200

\hsize=14.1cm
\vsize=19.5cm
\parindent=0cm   \parskip=0pt
\pageno=1

\def\ind{\hskip 1cm\relax}


\hoffset=15mm    
\voffset=1cm    


\ifnum\mag=\magstep1
\hoffset=-0.2cm   
\voffset=-.5cm   
\fi


\pretolerance=500 \tolerance=1000  \brokenpenalty=5000

\catcode`\@=11

\font\eightrm=cmr8         \font\eighti=cmmi8
\font\eightsy=cmsy8        \font\eightbf=cmbx8
\font\eighttt=cmtt8        \font\eightit=cmti8
\font\eightsl=cmsl8        \font\sixrm=cmr6
\font\sixi=cmmi6           \font\sixsy=cmsy6
\font\sixbf=cmbx6


\font\tengoth=eufm10       \font\tenbboard=msbm10
\font\eightgoth=eufm8      \font\eightbboard=msbm8
\font\sevengoth=eufm7      \font\sevenbboard=msbm7
\font\sixgoth=eufm6        \font\fivegoth=eufm5

\font\tencyr=wncyr10       
\font\eightcyr=wncyr8      
\font\sevencyr=wncyr7      
\font\sixcyr=wncyr6


\skewchar\eighti='177 \skewchar\sixi='177
\skewchar\eightsy='60 \skewchar\sixsy='60


\newfam\gothfam           \newfam\bboardfam
\newfam\cyrfam

\def\tenpoint{%
  \textfont0=\tenrm \scriptfont0=\sevenrm \scriptscriptfont0=\fiverm
  \def\rm{\fam\z@\tenrm}%
  \textfont1=\teni  \scriptfont1=\seveni  \scriptscriptfont1=\fivei
  \def\oldstyle{\fam\@ne\teni}\let\old=\oldstyle
  \textfont2=\tensy \scriptfont2=\sevensy \scriptscriptfont2=\fivesy
  \textfont\gothfam=\tengoth \scriptfont\gothfam=\sevengoth
  \scriptscriptfont\gothfam=\fivegoth
  \def\goth{\fam\gothfam\tengoth}%
  \textfont\bboardfam=\tenbboard \scriptfont\bboardfam=\sevenbboard
  \scriptscriptfont\bboardfam=\sevenbboard
  \def\bb{\fam\bboardfam\tenbboard}%
 \textfont\cyrfam=\tencyr \scriptfont\cyrfam=\sevencyr
  \scriptscriptfont\cyrfam=\sixcyr
  \def\cyr{\fam\cyrfam\tencyr}%
  \textfont\itfam=\tenit
  \def\it{\fam\itfam\tenit}%
  \textfont\slfam=\tensl
  \def\sl{\fam\slfam\tensl}%
  \textfont\bffam=\tenbf \scriptfont\bffam=\sevenbf
  \scriptscriptfont\bffam=\fivebf
  \def\bf{\fam\bffam\tenbf}%
  \textfont\ttfam=\tentt
  \def\tt{\fam\ttfam\tentt}%
  \abovedisplayskip=12pt plus 3pt minus 9pt
  \belowdisplayskip=\abovedisplayskip
  \abovedisplayshortskip=0pt plus 3pt
  \belowdisplayshortskip=4pt plus 3pt
  \smallskipamount=3pt plus 1pt minus 1pt
  \medskipamount=6pt plus 2pt minus 2pt
  \bigskipamount=12pt plus 4pt minus 4pt
  \normalbaselineskip=12pt
  \setbox\strutbox=\hbox{\vrule height8.5pt depth3.5pt width0pt}%
  \let\bigf@nt=\tenrm       \let\smallf@nt=\sevenrm
  \normalbaselines\rm}

\def\eightpoint{%
  \textfont0=\eightrm \scriptfont0=\sixrm \scriptscriptfont0=\fiverm
  \def\rm{\fam\z@\eightrm}%
  \textfont1=\eighti  \scriptfont1=\sixi  \scriptscriptfont1=\fivei
  \def\oldstyle{\fam\@ne\eighti}\let\old=\oldstyle
  \textfont2=\eightsy \scriptfont2=\sixsy \scriptscriptfont2=\fivesy
  \textfont\gothfam=\eightgoth \scriptfont\gothfam=\sixgoth
  \scriptscriptfont\gothfam=\fivegoth
  \def\goth{\fam\gothfam\eightgoth}%
  \textfont\cyrfam=\eightcyr \scriptfont\cyrfam=\sixcyr
  \scriptscriptfont\cyrfam=\sixcyr
  \def\cyr{\fam\cyrfam\eightcyr}%
  \textfont\bboardfam=\eightbboard \scriptfont\bboardfam=\sevenbboard
  \scriptscriptfont\bboardfam=\sevenbboard
  \def\bb{\fam\bboardfam}%
  \textfont\itfam=\eightit
  \def\it{\fam\itfam\eightit}%
  \textfont\slfam=\eightsl
  \def\sl{\fam\slfam\eightsl}%
  \textfont\bffam=\eightbf \scriptfont\bffam=\sixbf
  \scriptscriptfont\bffam=\fivebf
  \def\bf{\fam\bffam\eightbf}%
  \textfont\ttfam=\eighttt
  \def\tt{\fam\ttfam\eighttt}%
  \abovedisplayskip=9pt plus 3pt minus 9pt
  \belowdisplayskip=\abovedisplayskip
  \abovedisplayshortskip=0pt plus 3pt
  \belowdisplayshortskip=3pt plus 3pt
  \smallskipamount=2pt plus 1pt minus 1pt
  \medskipamount=4pt plus 2pt minus 1pt
  \bigskipamount=9pt plus 3pt minus 3pt
  \normalbaselineskip=9pt
  \setbox\strutbox=\hbox{\vrule height7pt depth2pt width0pt}%
  \let\bigf@nt=\eightrm     \let\smallf@nt=\sixrm
  \normalbaselines\rm}

\tenpoint


\def\pc#1{\bigf@nt#1\smallf@nt}         \def\pd#1 {{\pc#1} }






\def\^#1{\if#1i{\accent"5E\i}\else{\accent"5E #1}\fi}
\def\"#1{\if#1i{\accent"7F\i}\else{\accent"7F #1}\fi}



\newtoks\auteurcourant      \auteurcourant={\hfil}
\newtoks\titrecourant       \titrecourant={\hfil}

\newtoks\hautpagetitre      \hautpagetitre={\hfil}
\newtoks\baspagetitre       \baspagetitre={\hfil}

\newtoks\hautpagegauche
\hautpagegauche={\eightpoint\rlap{\folio}\hfil\the\auteurcourant\hfil}
\newtoks\hautpagedroite
\hautpagedroite={\eightpoint\hfil\the\titrecourant\hfil\llap{\folio}}

\newtoks\baspagegauche      \baspagegauche={\hfil}
\newtoks\baspagedroite      \baspagedroite={\hfil}

\newif\ifpagetitre          \pagetitretrue

a marche. Alors ...


\footline={\ifpagetitre\the\baspagetitre\else
\ifodd\pageno\the\baspagedroite\else\the\baspagegauche\fi\fi
\global\pagetitrefalse}


\def\raggedbottom{\topskip 10pt plus 36pt\r@ggedbottomtrue}



\def\pointir{\unskip . --- \ignorespaces}


\def\Bigbreak{\vskip-\lastskip\bigbreak}
\def\Medbreak{\vskip-\lastskip\medbreak}


\def\ctexte#1\endctexte{%
  \hbox{$\vcenter{\halign{\hfill##\hfill\crcr#1\crcr}}$}}


\long\def\ctitre#1\endctitre{%
    \ifdim\lastskip<24pt\vskip-\lastskip\bigbreak\bigbreak\fi
  		\vbox{\parindent=0pt\leftskip=0pt plus 1fill
          \rightskip=\leftskip
          \parfillskip=0pt\bf#1\par}
    \bigskip\nobreak}

\long\def\section#1\endsection{%
\vskip 0pt plus 3\normalbaselineskip
\penalty-250
\vskip 0pt plus -3\normalbaselineskip
\Bigbreak
\message{[section \string: #1]}{\bf#1\unskip}\pointir}

\long\def\sectiona#1\endsection{%
\vskip 0pt plus 3\normalbaselineskip
\penalty-250
\vskip 0pt plus -3\normalbaselineskip
\Bigbreak
\message{[sectiona \string: #1]}%
{\bf#1}\medskip\nobreak}

\long\def\subsection#1\endsubsection{%
\Medbreak
{\it#1\unskip}\pointir}

\long\def\subsectiona#1\endsubsection{%
\Medbreak
{\it#1}\par\nobreak}

\def\rem#1\endrem{%
\Medbreak
{\it#1\unskip} : }

\def\remp#1\endrem{%
\Medbreak
{\pc #1\unskip}\pointir}

\def\rema#1\endrem{%
\Medbreak
{\it #1}\par\nobreak}

\def\newparwithcolon#1\endnewparwithcolon{
\Medbreak
{#1\unskip} : }

\def\newparwithpointir#1\endnewparwithpointir{
\Medbreak
{#1\unskip}\pointir}

\def\newpara#1\endnewpar{
\Medbreak
{#1\unskip}\smallskip\nobreak}

\let\+=\tabalign

\def\signature#1\endsignature{\vskip 15mm minus 5mm\rightline{\vtop{#1}}}

\mathcode`A="7041 \mathcode`B="7042 \mathcode`C="7043 \mathcode`D="7044
\mathcode`E="7045 \mathcode`F="7046 \mathcode`G="7047 \mathcode`H="7048
\mathcode`I="7049 \mathcode`J="704A \mathcode`K="704B \mathcode`L="704C
\mathcode`M="704D \mathcode`N="704E \mathcode`O="704F \mathcode`P="7050
\mathcode`Q="7051 \mathcode`R="7052 \mathcode`S="7053 \mathcode`T="7054
\mathcode`U="7055 \mathcode`V="7056 \mathcode`W="7057 \mathcode`X="7058
\mathcode`Y="7059 \mathcode`Z="705A

\def\spacedmath#1{\def\packedmath##1${\bgroup\mathsurround=0pt ##1\egroup$}%
\mathsurround#1 \everymath={\packedmath}\everydisplay={\mathsurround=0pt }}

\def\nospacedmath{\mathsurround=0pt \everymath={}\everydisplay={} }


\long\def\th#1 #2\enonce#3\endth{%
   \Medbreak
   {\pc#1} {#2\unskip}\pointir{\it #3}\medskip}

\long\def\tha#1 #2\enonce#3\endth{%
   \Medbreak
   {\pc#1} {#2\unskip}\par\nobreak{\it #3}\medskip}


\long\def\Th#1 #2 #3\enonce#4\endth{%
   \Medbreak
   #1 {\pc#2} {#3\unskip}\pointir{\it #4}\medskip}

\long\def\Tha#1 #2 #3\enonce#4\endth{%
   \Medbreak
   #1 {\pc#2} #3\par\nobreak{\it #4}\medskip}


\def\decale#1{\smallbreak\hskip 28pt\llap{#1}\kern 5pt}
\def\decaledecale#1{\smallbreak\hskip 34pt\llap{#1}\kern 5pt}
\def\puce{\smallbreak\hskip 6pt{$\scriptstyle\bullet$}\kern 5pt}



\def\displaylinesno#1{\displ@y\halign{
\hbox to\displaywidth{$\@lign\hfil\displaystyle##\hfil$}&
\llap{$##$}\crcr#1\crcr}}


\def\ldisplaylinesno#1{\displ@y\halign{
\hbox to\displaywidth{$\@lign\hfil\displaystyle##\hfil$}&
\kern-\displaywidth\rlap{$##$}\tabskip\displaywidth\crcr#1\crcr}}


\def\eqalign#1{\null\,\vcenter{\openup\jot\m@th\ialign{
\strut\hfil$\displaystyle{##}$&$\displaystyle{{}##}$\hfil
&&\quad\strut\hfil$\displaystyle{##}$&$\displaystyle{{}##}$\hfil
\crcr#1\crcr}}\,}


\def\system#1{\left\{\null\,\vcenter{\openup1\jot\m@th
\ialign{\strut$##$&\hfil$##$&$##$\hfil&&
        \enskip$##$\enskip&\hfil$##$&$##$\hfil\crcr#1\crcr}}\right.}


\let\@ldmessage=\message

\def\message#1{{\def\pc{\string\pc\space}%
                \def\'{\string'}\def\`{\string`}%
                \def\^{\string^}\def\"{\string"}%
                \@ldmessage{#1}}}

\def\diagram#1{\def\normalbaselines{\baselineskip=0pt
\lineskip=10pt\lineskiplimit=1pt}   \matrix{#1}}



\def\up#1{\raise 1ex\hbox{\smallf@nt#1}}


\def\qed{\raise -2pt\hbox{\vrule\vbox to 10pt{\hrule width 4pt
                 \vfill\hrule}\vrule}}

\def\cqfd{\unskip\penalty 500\quad\vrule height 4pt depth 0pt width
4pt\medbreak}

\def\virg{\raise .4ex\hbox{,}}   


\def\build#1_#2^#3{\mathrel{
\mathop{\kern 0pt#1}\limits_{#2}^{#3}}}


\def\boxit#1#2{%
\setbox1=\hbox{\kern#1{#2}\kern#1}%
\dimen1=\ht1 \advance\dimen1 by #1 \dimen2=\dp1 \advance\dimen2 by #1
\setbox1=\hbox{\vrule height\dimen1 depth\dimen2\box1\vrule}%
\setbox1=\vbox{\hrule\box1\hrule}%
\advance\dimen1 by .4pt \ht1=\dimen1
\advance\dimen2 by .4pt \dp1=\dimen2  \box1\relax}

\def\date{\the\day\ \ifcase\month\or janvier\or f\'evrier\or mars\or
avril\or mai\or juin\or juillet\or ao\^ut\or septembre\or octobre\or
novembre\or d\'ecembre\fi \ {\old \the\year}}

\def\dateam{\ifcase\month\or January\or February\or March\or
April\or May\or June\or July\or August\or September\or October\or
November\or December\fi \ \the\day ,\ \the\year}

\def\crog{{\vrule height 2.57mm depth 0.85mm width 0.3mm}\kern -0.36mm
[}

\def\crod{]\kern -0.4mm{\vrule height 2.57mm depth 0.85mm
width 0.3 mm}}

\def\rond{\kern 1pt{\scriptstyle\circ}\kern 1pt}

\def\diagram#1{\def\normalbaselines{\baselineskip=0pt
\lineskip=10pt\lineskiplimit=1pt}   \matrix{#1}}

\def\hfl#1#2{\nospacedmath\smash{\mathop{\hbox to
12mm{\rightarrowfill}}\limits^{\scriptstyle#1}_{\scriptstyle#2}}}

\def\phfl#1#2{\nospacedmath\smash{\mathop{\hbox to
8mm{\rightarrowfill}}\limits^{\scriptstyle#1}_{\scriptstyle#2}}}

\def\vfl#1#2{\llap{$\scriptstyle#1$}\left\downarrow\vbox to
6mm{}\right.\rlap{$\scriptstyle#2$}}

\def\av{abelian variety}

\def\ppav{principally polarized abelian variety}
\def\ppavs {principally polarized abelian varieties}
\def\pa{\S\kern.15em}
\def\ra{\rightarrow}
\def\lra{\longrightarrow}
\def\llra{\nospacedmath\hbox to 10mm{\rightarrowfill}}
\def\lllra{\nospacedmath\hbox to 15mm{\rightarrowfill}}
\def\saut{\vskip 5mm plus 1mm minus 2mm}
\def\Z{\hbox{\bf Z}}
\def\P{\hbox{\bf P}}
\def\Q{\hbox{\bf Q}}
\def\C{\hbox{\bf C}}

\def\isom{\simeq}

\def\Supp{\mathop{\rm Supp}\nolimits}

\def\Pic{\mathop{\rm Pic}\nolimits}

\def\Id{\hbox{\rm Id}}
\def\rank{\mathop{\rm rank}\nolimits}

\def\ie{\hbox{i.e.}}
\def\resp{\hbox{resp.}}
\def\page{\hbox{p.}}

\def\A#1{{\cal A}_{#1}}

\def\J#1{{\cal J}_{#1}}
\def\CT{{\cal CT}_5}
\def\JT{(JC,\theta)}
\def\AT{(A,\theta)}

\def\cc#1{\hfill\kern .7em#1\kern .7em\hfill}

\def\dra{\ra\kern -3mm\ra}
\def\ldra{\lra\kern -3mm\ra}

\catcode`\@=12

\showboxbreadth=-1  \showboxdepth=-1\overfullrule=0pt
\baselineskip=14pt
\spacedmath{1.7pt}
\baspagegauche={\centerline{\tenbf\folio}}
\baspagedroite={\centerline{\tenbf\folio}}
\hautpagegauche={\hfil}
\hautpagedroite={\hfil}
\font\eightrm=cmr10 at 8pt
\font\pc=cmcsc10 \rm
\parskip=1.7mm

\def\CG#1{{\cal C}_{#1}}

\saut
\ctitre
{\bf MINIMAL COHOMOLOGY CLASSES AND JACOBIANS}
\endctitre
\medskip
\centerline{{\pc Olivier Debarre}
\footnote{(*)}{\rm Partially supported by
N.S.F. Grant DMS 92-03919 and the European Science
Project ``Geometry of Algebraic Varieties", Contract no. SCI-0398-C (A).}}

\vskip 1cm

\ind The main purpose of this article is to describe all effective algebraic
cycles
with minimal cohomology class in the Jacobian of a complex curve. More
precisely, let $\AT$
be a complex \ppav\ of dimension $g$. For $0\le d\le g$, the cohomology class
$\theta_d=\theta^d/d!$ is {\it minimal\/}, \ie\ non-divisible, in $H^{2d}(A,\Z
)$. When
$\AT$ is isomorphic to the Jacobian $\JT$ of a curve $C$ of genus $g$,
the image of the symmetric product
$C^{(g-d)}$ by any Abel-Jacobi map is a subvariety $W_{g-d}(C)$ of $JC$ with
class $\theta_d$ ([ACGH],
\page\ 25).

\ind Our main result (theorem 5.1) implies that {\it any
effective algebraic cycle in $JC$ with class $\theta_d$ is a translate of
either $W_{g-d}(C)$ or
$-W_{g-d}(C)$.}

\ind For $1<d<g$, I know of only one other family of \ppavs\ with an effective
algebraic
cycle with minimal class: in the $5$--dimensional intermediate
Jacobian $JT$ of a cubic threefold $T$ in $\P ^4$, the image by any Abel-Jacobi
map of the Fano
surface of lines contained in $T$ is a surface with class $\theta_3$ ([CG],
[B]).

\ind I suspect that these should be the only examples of effective algebraic
cycles with
class $\theta_d$ on abelian varieties of dimension $g$, when $1<d<g$. This
holds for any $g$ and
$d=g-1$ by Matsusaka's criterion ([M]), and for $g=4$ and $d=2$ by a result of
Ran ([R1]).

\ind In a second part, the main theorem is used to prove a weak version of this
conjecture: {\it for
$1<d<g$, the Jacobian locus (\resp\ the locus of  intermediate Jacobians of
cubic
threefolds) is an irreducible component of the set of \ppavs of dimension $g$
for which
$\theta_d$ (\resp\ $\theta_3$) is the class of an effective algebraic cycle}.
This result
was first proved by Barton and Clemens ([BC]) for $g=4$ and $d=2$; a slighty
weaker form of
it appeared later in [R1] (corollary III.1), but the proof is incomplete.

\ind Throughout this article, we will be working over the field of complex
numbers. A variety
will be a reduced projective (complex) scheme of finite type.

\ind The author would like to thank the M.S.R.I., where this research was done,
for
its hospitality and support.
\bigskip

\centerline{{\pc I. Subvarieties With Minimal Classes In Jacobians}}

\ind Let $\JT$ be the Jacobian of
a smooth curve $C$ of genus $g$. The aim of this part is to prove that any
effective
algebraic cycle in $JC$ with class $\theta_d$ is a translate of either
$W_{g-d}(C)$ or
$-W_{g-d}(C)$.

\ind We begin with a few preliminaries.

{\bf 1. Non-degenerate subvarieties}

\ind This notion was introduced by Ran ([R1], [R2]). Let $A$ be an
\av\ of dimension $g$ and let $W$ be a subvariety of $A$ of pure
dimension $d$. Let $W_{\rm reg}$ be the smooth part of $W$. We say that $W$ is
{\it non-degenerate} if the restriction map $H^0(A,\Omega^d_A)\lra H^0(W_{\rm
reg},\Omega^d_{W_{\rm reg}})$ is injective. By [R1], lemma II.1, this is
equivalent to each
one of the following properties:

\ind $\bullet$
the cup-product map $\cdot [W]:H^{d,0}(A)\lra H^{g,g-d}(A)$ is injective,
where\break $[W]\in
H^{g-d,g-d}(A)$ is the cohomology class of $W$,

\ind $\bullet$
the contraction map $\cdot \{ W\} :H^{g,d}(A)\lra H^{g-d,0}(A)$ is injective,
where\break $\{ W\} \in
H_{d,d}(A)$ is the homology class of $W$.

\ind This allows the extension of the definition to effective algebraic cycles.
Note that
on a \ppav ,  an effective algebraic cycle with class a multiple of a minimal
class is
non-degenerate.

{\bf 2. The property $({\cal P})$}

\ind   Let $V$ and $W$ be two irreducible subvarieties of an abelian variety
$A$ and let
$f:V\times W\ra A$ be the addition map. We will say that $V$ {\it has property
$({\cal P})$
with respect to $W$} if, for $v$ generic in $V$, the only irreducible
subvariety of
$f^{-1}(v+W)$ which dominates both $v+W$ via $f$, and $W$ via the second
projection, is $\{
v\}\times W$.

(2.1) Note that this implies that $\{
v\}\times W$ is a component of $f^{-1}(v+W)$, hence that the latter has the
same dimension
as $W$ at some point, hence that $f$ is generically finite onto its image; in
particular, $\dim
(V)+\dim (W)\le \dim (A)$. Note also that if $V$ has property $({\cal P})$ with
respect to
$W$, then  any translate of $V$ has property $({\cal P})$ with respect to any
translate of $W$.

{\pc Example} 2.2. If $C$ is a curve of genus $g$, it is not
difficult to check that $W_d(C)$ has property $({\cal P})$ with respect to
$W_e(C)$
whenever $d+e\le g$.

{\pc Lemma} 2.3. -- {\it Let $V$ and $W$ be two irreducible subvarieties of an
abelian
variety $A$ and let $g:W\times W\ra A$ be the subtraction map.
Then, $V$ has property $({\cal P})$ with respect to $W$ if and only if, for $v$
generic in $V$,
the only irreducible subvariety of $g^{-1}(V-v)$ which dominates both factors
$W$, is the
diagonal.}

{\bf Proof.} Let $f:V\times W\ra
A$ be the addition map. Let $v$ be generic in $V$ and consider the automorphism
$h$ of
$A\times A$ defined by $h(x,y)=(v+x-y,y)$. One checks that  $h\bigl(
g^{-1}(V-v)\bigr)=f^{-1}(v+W)$. The proposition then follows from the
fact that, if $Z$ is a variety
contained in $g^{-1}(V-v)$, then $f\bigl( h(Z)\bigr) = v+ p_1(Z)$ and  $p_2
\bigl(
h(Z)\bigr) = p_2(Z)$.\cqfd

{\pc Proposition} 2.4. -- {\it Let $V$ and $W$ be two irreducible subvarieties
of an abelian
variety $A$. Assume that $V$ has property $({\cal P})$ with respect to $W$.
Then:

\ind $\bullet$ the variety $(-V)$ has property $({\cal P})$ with respect to
$W$,

\ind $\bullet$ there exists a dense open set $\Omega$ in $V$ such that, if $U$
is an
irreducible subvariety of $V$ which meets $\Omega$, then $U$ has property
$({\cal P})$ with
respect to $W$.}

{\bf Proof.} The first  point follows from
lemma 2.3. Let $g:W\times W\ra A$ be the difference map. Let $\Omega$ be the
open set
 of points $v$ in $V$ such that the only irreducible subvariety of
$g^{-1}(V-v)$ which dominates both factors is the diagonal.  By the lemma,
$\Omega$ is dense in
$V$. Take $u$ in $\Omega\cap U$; since $U-u\subset V-u$, any irreducible
subvariety of
$g^{-1}(U-u)$ which dominates both factors is also a subvariety of
$g^{-1}(V-u)$, hence is
equal to the diagonal. By lemma 2.3, this proves the second point.\cqfd

\ind Although the following result can be checked directly, it is an easy
consequence of
lemma 2.3.

{\pc Proposition} 2.5. -- {\it Let $C$ be a curve of genus $g$ and let $V$ be
an
irreducible subvariety of $JC$. Assume that $V$ has property $({\cal P})$ with
respect to some
$W_d(C)$. Then  $V$ also has property $({\cal P})$ with respect to any
$W_e(C)$ for $0\le e\le d$.}

{\bf Proof.} Let $g:W_d(C)\times W_d(C)\ra JC$ and $h:W_e(C)\times W_e(C)\ra
JC$ be the
difference maps. Let $v$ be generic in $V$, and let $Z$ be an irreducible
subvariety of $h^{-1}(v-V)$ which
dominates both factors. Then $\{\ z+(w,w)\bigm| z\in Z\ ,\ w\in W_{d-e}(C)\ \}$
is
an irreducible subvariety of $g^{-1}(v-V)$ which dominates both factors. By
lemma 2.3, it is the
diagonal of $W_d(C)\times W_d(C)$, hence $Z$ is the diagonal of $W_e(C)\times
W_e(C)$. By lemma
2.3 again, this proves that $V$  has property $({\cal P})$ with respect to
$W_e(C)$.\cqfd

{\bf 3. Ran's theorem}

\ind Let $V$ and $W$ be two subvarieties of
 $A$, of respective pure codimensions $d$ and $g-d$. Assume that $s=V\cdot
W>0$. The
addition map $V\times W\ra A$ is then surjective, hence generically \'etale. It
follows that
for $x$ generic in $A$, the varieties $V$ and $x-W$ meet transversally at
distinct
smooth points $v_1,\ldots ,v_s$. Let $P_i: T_0A\ra T_0A$ be the projector with
image $T_{v_i}V$ and
kernel $T_{x-v_i}W$. Let
$ c(V,W)$ be the endomorphism $\sum_{i=1}^s \bigwedge ^{g-d} P_i$ of $\bigwedge
^{g-d} T_0A$;
Ran proves ([R1], theorem 2) that the transposed endomorphism ${}^tc(V,W)$
 of $\bigwedge ^{g-d} T_0^*A\isom H^{g-d,0}(A)$ is equal to the composition:
$$H^{g-d,0}(A)\buildrel{\cdot [V]}\over{\llra}H^{g,d}(A)\buildrel{\cdot
\{ W\} }\over{\llra}H^{g-d,0}(A)\ .$$
In particular, when $V$ and $W$ are non-degenerate, $c(V,W)$ is an
automorphism.

\ind The following
result of Ran on non-degenerate subvarieties with minimal intersection number
motivates
the somewhat abstruse definition of ``property $({\cal P})$''.

{\pc Theorem} 3.1. (Ran) -- {\it If $V$ and $W$ are non-degenerate subvarieties
of a
$g$--dimensional \ppav\ $\AT$, of respective pure codimensions $d$ and $g-d$,
then $V\cdot W\ge
{g\choose d}$. If moreover $W$ is irreducible and $V\cdot W={g\choose d}$, then
$V$ is
irreducible and has property $({\cal P})$ with respect to $W$.}

{\bf Proof.} Since $c(V,W)$ is the sum of $V\cdot W$ projectors of rank $1$, we
have:$${g\choose d}=\rank \bigl( c(V,W)\bigr)\le V\cdot W\ .$$
If there is equality, it follows from (7.4) in [R2] that $W$ meets at most one
component of $V$. But by [R1], corollary II.6, $W$ meets each component of $V$.
This
implies that $V$ is irreducible. To prove property $({\cal P})$, we keep the
same notation
as above and proceed as in the proof of [R1], theorem 5. If $\alpha_i$ spans
the line
$\bigwedge ^{g-d} T_{v_i}V$ in $\bigwedge ^{g-d} T_0A$, then $\{
\alpha_1,\ldots
,\alpha_s\}$ (with $s={g\choose d}$) is a basis for $\bigwedge ^{g-d} T_0A$.
The choice of
an identification $\bigwedge ^n T_0A\isom\C$ induces an isomorphism $ \bigwedge
^d
T_0A\isom\bigwedge ^{g-d} T_0^*A$. Let $\beta_i$ be an element
 of $\bigwedge ^d T_{x-v_i}W$ such that $\beta_i(\alpha_i)=1$. Then
$c(V,W)= \sum_{i=1}^s \alpha_i\otimes\beta_i$. In particular, for $i\ne j$:
$$\beta_i\bigl( c(V,W)^{-1}(\alpha_j)\bigr)=0\ .$$

\ind Fix $v=v_1$ in $V$ and let $x$
vary in $v+W$, so that $v\in V\cap (x-W)$. We get: $$\bigwedge ^d
T_{x-v_i}W\wedge
c(V,W)^{-1}\bigl(\bigwedge ^{g-d} T_vV\bigr)=0$$ for $i>1$. Since $W$ is
non-degenerate, the points
$x-v_1,\ldots ,x-v_s$ must therefore describe a proper subvariety of $W$. This
proves the last part
of the theorem, since $p_2f^{-1}(x)=W\cap (x-V)=\{ x-v,x-v_1,\ldots
,x-v_s\}.$\cqfd

\ind For example, if $C$ is a curve of genus $g$, the theorem implies that
$W_d(C)$ has
 property $({\cal P})$ with respect to $W_{g-d}(C)$ for $d\le g$. Together with
proposition 2.5, this proves the claim of example 2.2.

{\bf 4. An auxiliary result}

\ind Let $C$ be a smooth
curve of genus $g$ and let $n>1$. We define a subset ${\cal T}_n$ of $C^{(n)}$
as
follows: consider all surjective morphisms $C\ra C'$ of degree $r>1$, where
$C'$ is a smooth
irrational curve. Such a morphism induces a map $\psi:C'\ra C^{(r)}$; let
${\cal T}_n$
be the union of all $\psi(C')+C^{(n-r)}$, for $r\le n$, obtained in this way,
of the inverse image in $C^{(n)}$
of $W^1_n(C)$ and of
the diagonal $2C+C^{(n-2)}$. For $n\le g$, it is a proper closed subvariety of
$C^{(n)}$.

{\pc Proposition} 4.1. -- {\it Let $\JT$ be the Jacobian of a smooth curve $C$
of genus
$g\ge 2$ and let $Z$ be a subvariety of $C^{(n)}$ of pure codimension $m$.
Then, for any
$E$ in $C^{(n-m)}$, the variety $Z$ meets $E+C^{(m)}$. Moreover, if there
exists $E$ in
$C^{(n-m)}$ such that $Z\cdot (E+C^{(m)})=1$, then either $Z$ is
contained in ${\cal T}_n$, or
there exist points $c_1,\ldots,c_m$ of $C$ such that
$Z=c_1+\ldots+c_m+C^{(n-m)}$.}

{\bf Proof.} The first part follows from the fact that the cohomology class of
$E+C^{(m)}$
is the $(n-m)$--fold self-intersection of the cohomology class of the ample
divisor
$C^{(n-1)}$ in $C^{(n)}$ ([ACGH], p\page\ 309, 310).

\ind Assume now that $Z\cdot (E+C^{(m)})=1$. We first do the case
when $Z$ is a curve. We may assume that $Z$ is not contained in $x+C^{(n-1)}$
for any $x$ in $C$.
Since $Z\cdot (x+C^{(n-1)})=1$ for any $x$ in $C$, the curve $Z$ is {\it
smooth} and there exists a
morphism $\tau:C\ra C^{(n-1)}$ such that $Z=\{ x+\tau (x)\bigm| x\in C\} $. Let
$\Gamma$ be the curve
$\{ (x,\tau (x))\bigm| x\in C\} $ in $C\times
C^{(n-1)}$. If the induced morphism $\phi:\Gamma\ra Z$ is not birational, $Z$
is
contained in ${\cal T}_n$.

\ind Since $Z$ is smooth, we are left with the case where $\phi$ is an
isomorphism.
If $\tau$ is constant, the proof is over. Assume therefore that $\tau(C)$ is a
curve. If $n=2$, we get
$Z=\{ x+\tau (x)\bigm| x\in C\} $, where $\tau$ is an involution of $C$. In
particular, $Z$
is contained in ${\cal T}_2$. We assume $n\ge 3$ and proceed by induction. Let
us show that
the curve $\tau(C)$ satisfies the same property as $Z$. Let $x$ be a point of
$C$ and assume
that $\tau(y)=x+D$ and $\tau(y')=x+D'$ are both in $x+C^{(n-1)}$ for some $y$
and $y'$ on $C$.
Then, $y+x+D$ and $y'+x+D'$ are on $Z$ hence the hypothesis on $Z$ implies
$y+D=y'+D'$. Therefore,
either $y=y'$ and $D=D'$, or there exists an element $E$ of $C^{(n-3)}$ such
that $D=y'+E$ and
$D'=y+E$. Then, both $(y,\tau(y))$ and $(y',\tau(y'))$ are sent by addition to
$x+y+y'+E$
hence, since $\phi$ is an isomorphism, we get again $y=y'$ and $D=D'$. It
follows that $\tau(C)$ and $x+C^{(n-2)}$ have a single common point and one
checks
that the intersection is transverse. We can therefore apply the induction
hypothesis to
$\tau(C)$: since we have assumed that $Z$ is not contained in any
$x+C^{(n-1)}$, the curve $\tau(C)$
 is in ${\cal T}_{n-1}$. It follows that $Z$ is contained in $C+{\cal
T}_{n-1}$,
hence in ${\cal T}_n$. This
finishes the proof of the lemma when $Z$ is a curve.

\ind We now do the general case, by induction on the dimension of $Z$. Let $Z$
of
dimension $n-m>1$ satisfy the hypothesis. For $x$ generic in $C$,
the inverse image of $Z$ by the map $C^{(n-1)}\ra C^{(n)}$ which sends $D$ to
$x+D$
satisfies the same property. Therefore, either $Z$ is contained in $C+{\cal
T}_{n-1}$, hence
in ${\cal T}_n$, or there exists a morphism $\tau:C\ra C^{(m)}$ such that:
$$Z=\{\ x+\tau(x)+C^{(n-m-1)}\bigm| x\in C\ \}\ .$$
But $\tau$ must then be
constant. This finishes the proof of the proposition.\cqfd

\medskip
{\bf 5. The main theorem}

\ind We now prove our main result.

{\pc Theorem} 5.1. -- {\it Let $\JT$ be the Jacobian of a smooth curve $C$ of
genus $g$ and
let $V$ be an effective non-degenerate algebraic $(g-d)$--cycle in $JC$ such
that $V\cdot
W_d(C)={g\choose d}$. Then $V$ is a translate of either $W_{g-d}(C)$ or
$-W_{g-d}(C)$.}

\ind  We refer to \S 1 for the definition of a non-degenerate cycle. Recall
that any effective
algebraic cycle with class $\theta_d$ satisfies the hypotheses of the theorem.

{\bf Proof.} By replacing any multiple component of $V$ by a sum of translates,
we may assume
that $V$ is a subvariety. Ran's theorem 3.1 implies that $V$ is irreducible and
has property
$({\cal P})$ with respect to $W_d(C)$. We first prove by  induction on $g-d$
that any
$(g-d)$--dimensional irreducible subvariety $V$ of $JC$ which has property
$({\cal P})$ with respect to $W_d(C)$ is a translate of some
$W_r(C)-W_{g-d-r}(C)$, with
$0\le r\le g-d$. This is obvious for $g-d=0$, hence we assume $g-d>0$.

\ind For any positive integer $e$, we define the
addition map $f^+_e:V\times W_e(C)\ra JC$ and the subtraction map
$f^-_e:V\times W_e(C)\ra JC$. When $e\le d$, proposition 2.5 implies that $V$
has property
$({\cal P})$ with respect to $W_e(C)$, and so does $-V$ (proposition
2.4); it follows from (2.1) that $f^+_e$ and
$f^-_e$ are finite onto their images.

{\pc Proposition} 5.2. -- {\it Let $V$ be an irreducible subvariety of $JC$ of
codimension $d<g$
having property $({\cal P})$ with respect to $W_d(C)$. Then there exists an
irreducible
subvariety $U$ of $JC$ such that either $V=U+C$ or $V=U-C$.}

{\bf Proof.} The maps $f^{\pm}_d$ are not birational ([R1], lemma II.15). Let
$n$ be the
smallest integer such that either $f^+_n$ or $f^-_n$ is not birational onto its
image. Since
both the hypotheses and conclusion of the lemma hold for $V$ if and only if
they do for $-V$
(proposition 2.4), we may assume that $f^+_n$ is not birational onto its image.
For $v$
generic in $V$, there exists an irreducible component $\Gamma$ of
$(f^+_n)^{-1}(v+W_n(C))$,
distinct from $\{ v\} \times W_n(C)$, which dominates $v+W_n(C)$. Since $V$ has
property
$({\cal P})$ with respect to $W_n(C)$ (proposition 2.5), the projection
$p_2(\Gamma)$ is a subscheme $G$ of $W_n(C)$ of dimension $m<n$. For $D$
generic in $G$,
there exists a subvariety $S_D$ of $V$ of dimension $n-m$ such that: $$
S_D+D\subset v+W_n(C)\ .$$
\ind Let $E$ be generic in $W_m(C)$. By proposition 4.1, the subvariety
$E+W_{n-m}(C)$ of
$W_n(C)$ meets $G$ at some point $D=E+E'$. For the same reason, the subvariety
$S_D+D-v$ of $W_n(C)$
meets $E'+W_m(C)$ at $r$ points $E'+E_1,\ldots, E'+E_r$ with $r>0$, \ie\ there
exist points
$v_1,\ldots,v_r$ of $S_D\subset V$ such that: $$
v_i+(E+E')-v=E'+E_i\ ,\qquad {\rm \ie}\qquad
v_i-E_i=v-E\ .$$
\ind Since $m<n$, the map $f^-_m:V\times W_m(C)\ra JC$ is
birational onto its image. Therefore, since $v$ and $E$ are generic, we must
have $E_i=E$ for all
$i$ and $v\in S_D$. In other words: $$
(S_D+D-v)\cap \bigl( E'+W_m(C)\bigr) =\{ D\}$$
as sets. This is actually a scheme-theoretic equality; if $\epsilon$ is a
tangent vector
such that:
{\nospacedmath
$$\displaylines{
\epsilon\in T_D(S_D+D-v)=T_vS_D\subset T_vV\cr
\epsilon\in T_D\bigl( E'+W_m(C)\bigr) =T_EW_m(C)\ ,\cr}$$}then
$(\epsilon,\epsilon)$ is in the
kernel of the differential of $f^-_m$ at $(v,E)$, hence is $0$.

\ind By proposition 4.1, $E'+W_m(C)$ meets each component of $S_D+D-v$, which
must therefore
be irreducible. Recall that since $f^+_n(\Gamma)=v+W_n(C)$, the subvarieties
$S_D+D-v$
cover $W_n(C)$ as $D$ varies in $G$. In particular, we may assume that
$S_D+D-v$ is not
contained in the image of ${\cal T}_n$ (defined in \S 4) in $W_n(C)$. It then
follows from
proposition 4.1 applied to the strict transform of $S_D+D-v$ in $C^{(n)}$ that
there exists
an effective divisor $E_D$ of degree $m$ such that:
$$
S_D+D-v=W_{n-m}(C)+E_D\ .$$
Moreover, since $v\in S_D$, the divisor $E'_D=D-E_D$ is effective. If $c$ is in
the support
of $E'_D$ and if $x$ is any point of $C$, the point $v-c+x$ is in $S_D$ hence
in $V$. Since $v$ is
generic, this implies that some
irreducible component $T$ of the scheme: $$\{\ ( v,c)\in V\times C\bigm|
v-c+C\subset V\ \}$$
dominates $V$. Note that $U=f^-_1(T)$ satisfies $U+C\subset V$. In particular,
$\dim U<$\break $\dim
V\le\dim T$. In general, if $f^-_1$ contracts a subvariety $F$ of $V\times C$
to a point, the
projection $F\ra C$ is an isomorphism. This implies that $U= f^-_1(T)$ has
dimension $\dim T -1$,
and that $V=p_1 (T) \subset U+ C$, hence $V= U+ C$.\cqfd

\ind Replacing $V$ with $-V$ if necessary, we may assume that $V=U+C$. The
following
proposition shows that $U$ satisfies the induction hypothesis.

{\pc Proposition} 5.3. -- {\it Let $V$ be an irreducible subvariety of $JC$ of
codimension $d<g$
 having property $({\cal P})$ with respect to $W_d(C)$. Assume that $V=U+C$.
Then $U$
has property $({\cal P})$ with respect to $W_{d+1}(C)$.}

{\bf Proof.} By proposition 2.4, $U$ has property $({\cal
P})$ with respect to $W_d(C)$, hence also with respect to $C$ (proposition
2.5).

\ind Let $u$ be generic in $U$ and let $\phi^+_{d+1}:U\times
W_{d+1}(C)\ra JC$ be the addition map. Let $\Gamma_0$ be an irreducible
subvariety of $U\times
W_{d+1}(C)$ such that $\phi^+_{d+1}(\Gamma_0)=u+W_{d+1}(C)$ and
$p_2(\Gamma_0)=W_{d+1}(C)$. Let $\Gamma$
be an irreducible subvariety of $U\times C^{(d+1)}$ which maps onto $\Gamma_0$
by the natural map:
$$\Pi_U=\Id _U\times \pi_{d+1}:U\times C^{(d+1)}\ra U\times
W_{d+1}(C)\ .$$We have $\phi^+_{d+1}\Pi_U (\Gamma)=u+W_{d+1}(C)$ and
$p_2(\Gamma)=C^{(d+1)}$. Moreover,
since $\phi^+_{d+1}\Pi_U$ is generically finite onto its image, the dimension
of $\Gamma$ is $d+1$. We
need to show that $\Gamma=\{ u\}\times C^{(d+1)}$.

\ind For any $c$ in $C$, each component of $\Gamma\cap
(\phi^+_{d+1}\Pi_U)^{-1}\bigl( u+c+W_d(C)\bigr)$ that dominates $u+c+W_d(C)$
has dimension $d$. Since
$p_2(\Gamma)=C^{(d+1)}$, at least one of these components, say $\Gamma_c$,
projects onto a divisor
$Z_c$ in $C^{(d+1)}$. We use the elementary:

{\pc Lemma} 5.4. -- {\it Let $Z$ be an irreducible divisor in $C^{(d+1)}$ and
let $g:C\times C^{(d)}\ra C^{(d+1)}$ be the natural map. Then, there exists
an irreducible divisor $Z'$ in
$C\times C^{(d)}$ such that $g(Z')=Z$ and $p_2(Z')=C^{(d)}$.}

{\bf Proof.} By proposition 4.1,
$g^*Z\cdot (C\times\{ D\} )=Z\cdot (C+D)$ is non-zero for any point $D$ in
$C^{(d)}$, hence
$p_2\bigl( g^{-1}(Z)\bigr) =C^{(d)}$. Some component $Z'$ of $g^{-1}(Z)$ must
satisfy the same
property. Its dimension is then $\ge d$, hence it must map onto $Z$ by the
finite map $g$.\cqfd

\ind Pick such a divisor $Z'_c$ for $Z_c$ and let $\Gamma'_c$ be a component of
$Z'_c\times
_{Z_c}\Gamma_c$ that projects onto $Z'_c$. The situation is the following:
$$\diagram{\Gamma'_c&\subset&U\times C\times C^{(d)}&\hfl{\Id _U\times
g}{}&U\times
C^{(d+1)}&\supset&\Gamma_c\cr
\vfl{}{}&&\vfl{p_{23}}{}&&\vfl{p_2}{}&&\vfl{}{}\cr
Z'_c&\ \subset&C\times C^{(d)}&\hfl{g}{}&C^{(d+1)}&\supset&Z_c\cr
\vfl{}{}&&\vfl{p_2}{}&&&&\cr
C^{(d)}&\ \ =&C^{(d)}&&&&\cr}
$$
Since $\Id _U\times g$ is finite, it maps $\Gamma'_c$ onto $\Gamma_c$. Set
$\phi^+_{1,d}=\phi_1^+\times\Id _{C^{(d)}}$ and consider: $$U\times C\times
C^{(d)}\buildrel
\phi^+_{1,d}\over {\llra}V\times C^{(d)}\buildrel \Pi_V\over {\llra} V\times
W_d(C)\buildrel
\phi^+_d\over {\llra}JC\ . $$
Then:
{\nospacedmath
$$\displaylines{
\phi^+_d\Pi_V\phi^+_{1,d}(\Gamma'_c)=\phi^+_{d+1}\Pi_U (\Gamma_c)=u+c+W_d(C)\cr
p_2\Pi_V\phi^+_d(\Gamma'_c)=\pi_d p_3(\Gamma'_c)=W_d(C)\ .\cr}$$}Since $V$ has
property
$({\cal P})$ with respect to $W_d(C)$, this implies:
$$\Pi_V\phi^+_{1,d}(\Gamma'_c)=\{
u+c\}\times W_d(C)$$ for $c$ generic. Therefore:
 $$\Gamma_c=\{ u'\}\times (c'+C^{(d)})\qquad {\rm and}\qquad Z_c=c'+C^{(d)}\
,$$
for some points $u'$ of $U$ and $c'$ of $C$ such that $u'+c'=u+c$. The union as
$c$ varies in $C$ of
all $Z_c$'s must be $C^{(d+1)}$. It follows that $c'$ must describe the whole
of $C$. Recall that $U$
has property ${\cal P}$ with respect to $C$, hence that the only component of
$(\phi_1^+)^{-1}(u+C)$ which dominates both $u+C$ via $\phi_1^+$ and $C$ via
$p_2$ is $\{
u\}\times C$. This implies $\Gamma_c=\{ u\}\times (c+C^{(d)})$ for $c$ generic,
hence the
proposition.\cqfd

\ind It then follows from our induction hypothesis that $U$ is a
translate of some $W_{r-1}(C)-W_{g-d-r}(C)$, hence that
$V$ is a translate of $W_r(C)-W_{g-d-r}(C)$.

\ind When $C$ is hyperelliptic, $-C$ is a translate of $C$, hence
$-W_{g-d-r}(C)$ is a
translate of $W_{g-d-r}(C)$, and $V$ is a translate of $W_{g-d}(C)$. To
conclude the proof of
the theorem, we may therefore assume that $C$ is non-hyperelliptic.

{\pc Lemma} 5.5. -- {\it Take $r,d\ge 0$ with $d<g$ and $d+r\le g$, and assume
that $C$ is
non-hyperelliptic. Then:
$$\bigl( W_r(C)-W_{g-d-r}(C)\bigr)\cdot W_d(C)={g-d\choose r}{g\choose d}\ .$$}
{\bf Proof.} The subtraction map $W_r(C)\times W_{g-d-r}(C)\ra JC$ is
birational onto its
image. In fact, if $(D,E)$ and $(D',E')$ are two distinct elements of
$W_r(C)\times W_{g-d-r}(C)$ such that $D-E\equiv D'-E'$, then:

\ind $\bullet$ either $h^0(C,D+E')=1$, in which case $(D,E)$ belongs
to the divisor\break $\{ (D,E)\,\bigm|\,\Supp (D)\cap\Supp (E)\ne \emptyset
\}$,

\ind $\bullet$ or $h^0(C,D+E')>1$, in which case $(D,E)$ belongs to: $$\{\
(D,E)\,\bigm|\,\exists L\in W^1_{g-d}(C)\ \ \ D,E\le |L|\ \}\ .$$
This locus has dimension
$\max _{a>0}\bigl(\dim W^a_{g-d}(C)\, +2a\bigr)$, which is less than $g-d$ if
$C$ is
non-hyperelliptic (Martens' theorem, [ACGH], \page\ 191).

\ind If $\mathrel{*}$ denotes the Pontryagin product on $H^{\bullet}(JC,\Q )$,
it follows that
the cohomology class of $\bigl( W_r(C)-W_{g-d-r}(C)\bigr)$ is:
$$\theta_{g-r}\mathrel{*}\theta_{d+r}={g-d\choose r}\theta_d\ .$$
Its intersection number with $W_d(C)$ is therefore:
$${g-d\choose r}\theta_d\cdot\theta_{g-d}={g-d\choose r}{g\choose d}\ .$$
\ind This proves the lemma.\cqfd

\ind But this intersection number is equal to $V\cdot W_d(C)$, hence to
${g\choose
d}$. It follows that either $r=g-d$, or $r=0$, hence either $V$ is a translate
of
$W_{g-d}(C)$, or it is a translate of $-W_{g-d}(C)$.\cqfd

\vfill\supereject\null

\centerline{{\pc II. A Weak Characterization Of Jacobians}}

\medskip
{\bf 6. Preliminaries}

\ind For any positive integer $g$, let $\A g$ be
the moduli space of complex \ppavs of dimension $g$,
let $\J g$ be the closure in $\A g$ of the subvariety which corresponds
to Jacobians of smooth curves of genus $g$, and let $\CT$ be the closure in $\A
5$ of the
subvariety which corresponds to intermediate Jacobians of smooth cubic
threefolds. For
$0<d\le g$, the subset $\CG {g,d}$ of $\A g$ which corresponds to
\ppavs $\AT$ for which $\theta_d$ is the class of an effective algebraic cycle
(so that
$\CG {g,1} =\CG {g,g} = \A g$), is {\it closed} in $\A g$. Indeed, let ${\cal
X}\ra S$ be
a versal family of abelian varieties of dimension $g$ with a relatively
ample line bundle ${\cal L}$ on ${\cal X}$ which induces principal
polarizations on the
fibers, such that the classification morphism $S\ra \A g$ is a finite cover.
Fix an embedding of
${\cal X}$ in some projective space $\P ^N_S$, using for example the sections
of ${\cal
L}^{\otimes 3}$. Chow coordinates show that the family of effective cycles in
${\cal X}$ whose
fibers over $S$ have codimension $d$ and degree $\theta_d(3\theta)^{g-d}$, is
projective over
$S$. In particular, so is its closed subset which parametrizes effective cycles
with class
$\theta_d$ in the fibers. Consequently, its image in $S$ is also closed, hence
so is its image
$\CG {g,d}$ in $\A g$.

 \medskip
{\bf 7. Degenerations of abelian varieties}

\ind Let $\Delta=\{ t\in\C \bigm| |t|<1\}$. A degeneration of \ppavs of
dimension $g+1$ will
be a proper family ${\cal X}\ra\Delta$, with ${\cal X}$ smooth, whose fibers
over
$\Delta^*$ are \ppavs of dimension $g+1$ and whose central fiber $X$ is a
projective
variety whose normalization is a $\P ^1$--bundle $\P$ over a \ppav\
$(A,\theta)$ of dimension
$g$. More precisely, there exists an element $a$ of $A$ such that $\P =\P
\bigl({\cal
O}_A\oplus{\cal O}_A(\widehat a)\bigr)$, where $\widehat a$ is the image of $a$
by the
canonical isomorphism $\phi_{\theta}:A\ra\Pic ^0(A)$. The bundle $p:\P\ra A$
has two
disjoint sections $\P _0$ and $\P _{\infty}$, and $X$ is obtained from $\P$ by
identifying
any point $x$ on $\P _0$ with the point $x-a$ on $\P _{\infty}$. The singular
locus $X_s$
of $X$ will always be identified with $A$ via the isomorphism $X_s\ra\P _0\ra
A$. By [N],
theorem 16.1, we may
also assume the existence of a relatively ample line bundle ${\cal L}$ on
${\cal X}$ which
induces the principal polarization on the smooth fibers and which restricts
 on $X$ to an ample line bundle whose pull-back to $\P$ is ${\cal
O}_{\P}(p^*\Theta+\P _0)\simeq{\cal O}_{\P}(p^*\Theta_a+\P _{\infty}) $, for a
suitable
representative $\Theta$ of the polarization $\theta$. Another (more accessible)
reference for
the existence of ${\cal L}$ is [HW], proposition 4.1.3 (the proof is given in
the surface
case, but works in general).

\ind The projective varieties $X$ with their polarizations correspond to the
boundary
points of $\A {g+1}$ in a suitable partial compactification $\bar{{\cal
A}}_{g+1}$,
where they form a divisor $\partial\A {g+1}$. There is a
surjective map $q:\partial\A {g+1} \ra\A g$ which, in the above notation, sends
$X$ with
its polarization to $(A,\theta)$, and whose generic fiber is isomorphic to
$A/\pm 1$. For any
subvariety ${\cal B}$ of $\A {g+1}$, we will denote by $\partial{\cal B}$ the
intersection
of $\partial\A {g+1}$ with the closure of ${\cal B}$ in $\bar{{\cal A}}_{g+1}$.

{\pc Proposition} 7.1. -- {\it For $1<d<g+1$, the image by $q$ of the boundary
$\partial\CG
{g+1,d}$ is contained in $\CG{g,d}\cap\CG{g,d-1}$.}

{\bf Proof.} Keeping the same notation, we assume moreover that the smooth
fibers of ${\cal
X}\ra\Delta$ correspond to elements of $\CG {g+1,d}$. Then,
 there exists a subvariety ${\cal Z}$ of ${\cal X}$ such that
the cycles associated with the fibers of ${\cal Z}\ra\Delta$ over $\Delta^*$
are sums of
effective cycles, each with class $\theta_d$. Over the open interval $(0,1)$,
the scheme ${\cal
Z}$ contains a subscheme whose fibers have class $\theta_d$.  Its topological
closure
meets the central fiber $X$ in a union $Z=Z_1\cup\cdots\cup Z_s$ of irreducible
algebraic
subvarieties of codimension $d$. For $i=1,\ldots ,s$, let $Y_i$ be the closure
in $\P$ of
the inverse image of $Z_i$ by the bijection $\P -\P _{\infty}\ra X$. Then:
$$
\sum_{i=1}^s m_i [Y_i] =\bigl( p^*\theta +[\P _0]\bigr) ^d/d!
 = p^*\theta_d\ + \ [\P _0]\cdot p^*\theta_{d-1}\ , \leqno (7.2)$$

for some positive integers $m_1,\ldots ,m_s$. It follows that the image of
$\sum_{i=1}^s m_i Y_i$ in $A$ has class $\theta_{d-1}$, and that its
intersection with $\P _{\infty}$ has class $\theta_d$. This proves the
proposition.\cqfd

\medskip
{\bf 8. Proof of the theorem}

{\pc Theorem} 8.1. -- {\it For $1<d<g$, the locus $\J g$ is a component of $\CG
{g,d}$.
Moreover, $\CT$ is a component of $\CG {5,3}$ and is not contained
in $\CG {5,2}$.}

{\bf Proof.} The proof of the first part is by induction on $g$. Since this
property is empty for
$g=2$, we assume that it holds for some $g\ge 2$ and show that it then holds in
dimension $g+1$. Let
$1< d< g+1$ and let ${\cal F}$ be a component of $\CG {g+1,d}$ which contains
either $\J {g+1}$ or
$\CT$ (if $g=4$). It is known that $q(\partial\J {g+1}) = \J g$ and
$q(\partial\CT) = \J 4$. It
follows from proposition 7.1 and our induction hypothesis that $\J g$ is a
component of
$q(\partial{\cal F})$. With the notation of the proof of proposition 7.1, we
may therefore
assume that $(A,\theta)$ is the Jacobian $(JC,\theta)$ of a generic smooth
curve $C$ of
genus $g$. Then, $\partial\J g\cap q^{-1}\{ (JC,\theta)\}$ corresponds to the
image of the
surface $C-C$ in $JC/\pm 1$. For $g=4$, the curve $C$ has two pencils $g^1_3$
and $h^1_3$ of
degree $3$, and $\partial\CT\cap q^{-1}\{ (JC,\theta)\}$ corresponds to the
point $\pm
(g^1_3-h^1_3)$ of $JC/\pm 1$ ([C]).

\ind We will show that $\partial\CG {g+1,d}\cap q^{-1}\{ (JC,\theta)\}$ is
equal to
$(C-C)\cup \{ \pm
(g^1_3-h^1_3)\}$ for $g=4$ and $d=3$, and to $C-C$ otherwise.

\ind Keeping the notation of the proof of proposition 7.1, the
image in $JC$ of the effective cycle $Y=\sum_{i=1}^s m_i Y_i$ has class $\theta
_{d-1}$. By theorem 5.1, it is a translate of $\pm W_{g-d+1}(C)$ which,
possibly after
translation and action of $(-1)_X$, we may assume to be
$W_{g-d+1}(C)$. Consequently, after reindexing if necessary, we have $m_1=1$,
the morphism
$Y_1\ra p(Y_1)$ is birational, and $Y_i=p^{-1}\bigl( p(Y_i)\bigr)$ for $i>1$.
Identity (7.2)
and theorem 5.1 then imply that:

\ind {\bf Case 1:} either $Y=Y_1$ maps birationally onto $W_{g-d+1}(C)$ and
meets $\P
_0$ and $\P _{\infty}$ along subvarieties which are translates of $\pm
W_{g-d}(C)$,

\ind {\bf Case 2:} or $Y=Y_1+Y_2$, where $Y_1$ has class $[\P _0]\cdot
p^*\theta
_{d-1}$ and maps isomorphically onto $W_{g-d+1}(C)$, and $Y_2$ is the pull-back
by $p$ of
a translate $V$ of $\pm W_{g-d}(C)$. In this case, either $Y_1$ is contained
in $\P _0$ (which means that $Z$ has a component contained in the singular
locus $X_{s}$ of
$X$, which may well happen, contrary to what is asserted in [BC] on \page\ 64),
or it meets
neither $\P _0$ nor $\P _{\infty}$. In the latter case, the line bundle ${\cal
O}_{\P }(\P
_0 - \P _{\infty})=p^*{\cal O}_{JC}(\widehat a )$ is trivial on $Y_1$. Since
the restriction
$\Pic ^0(JC)\ra\Pic ^0\bigl( W_{g-d+1}(C)\bigr) $ is injective, this implies
that $a$
is $0$ hence is in $C-C$. We
may therefore assume that $Y_1$ is contained in $\P _0$.

\ind The following result makes the structure of the limit cycle more precise.

{\pc Lemma} 8.2. -- {\it Let $z$ be a point of $Z\cap X_{s}$. Then either $z$
belongs to a component of $Z$ contained in $X_s$, or each of the two branches
of
$X$ contains a local component of $Z$ at $z$.}

{\bf Proof.} Let $\varepsilon: {\widetilde {\cal X}}\ra {\cal X}$ be the
blow-up of the
smooth subscheme $X_s$. The strict transform ${\widetilde X}$ of $X$ can be
identified
with $\P$. Let ${\widetilde {\cal Z}}$ be the strict transform of ${\cal Z}$.
If
no component of $Z$ through $z$ is contained in $X_s$, the scheme $Z\cap
X_s={\cal Z}\cap
X_s$ has pure codimension $2$ in ${\cal Z}$ at $z$.
It follows that its ideal in ${\cal Z}$ is not invertible at $z$, hence that
the rational
curve $\varepsilon^{-1}(z)$ is contained in ${\widetilde {\cal Z}}$. This curve
meets
${\widetilde X}$ at two points (one on $\P _0$, the other on $\P _{\infty}$)
which
correspond to the two branches of $X$ at $z$. It follows that the intersection
of
${\widetilde {\cal Z}}$ with the Cartier divisor ${\widetilde X}$ is non-empty,
hence has
dimension $\dim {\widetilde {\cal Z}}-1=g+1-d$, at each of these two points.
Consequently,
the intersection of $Z$ with each of the two branches of $X$ has dimension
$g+1-d$ at $z$,
and this proves the lemma.\cqfd

\ind In case 1, this means that $V=Y\cap \P _0$ and $Y\cap \P _{\infty}$ must
be identified through the glueing process, hence that $Y\cap \P _{\infty}=V-a$.
In
particular, both $V$ and $V-a$ are contained in $W_{g-d+1}(C)$.

\ind In case 2, both $V$ and $V-a$ must be again contained in
$W_{g-d+1}(C)$.

\ind Recall that in both cases, $V$ is a translate of $\pm W_{g-d}(C)$.

{\pc Lemma} 8.3. -- {\it Let $C$ be a smooth curve of genus $g$ and let $a$ be
an
element of $JC$. Assume that some translate of $W_{m+1}(C)$ contains
$\varepsilon W_m(C)$
and $\varepsilon W_m(C)-a$, for some $\varepsilon=\pm 1$ and $0\le m\le g-2$.
Then either
$a\in C-C$, or $\varepsilon=-1$, $g=4$ and $a$ is the difference of two
$g^1_3$'s.}

{\bf Proof.} Let $L_a$ be the line bundle of degree $0$ on $C$ associated with
$a$. If
$\varepsilon=1$, there exists a line bundle $L$ on $C$ of degree $1$ such that:
$$h^0\bigl( C,L(c_1+\cdots +c_m)\bigr) >0\qquad{\rm and}\qquad
h^0\bigl( C,L\otimes L^{-1}_a(c_1+\cdots +c_m)\bigr) >0\ ,$$for any points
$c_1,\ldots,c_m$ of $C$. Riemann-Roch implies that both $L$ and $L\otimes
L^{-1}_a$ have a
non-zero-section, hence that $a\in C-C$.

\ind If $\varepsilon=-1$, we may assume that $C$ is non-hyperelliptic and
$0<m\le g-3$
(otherwise, $-W_{m+1}(C)$ is
 a translate of $W_{m+1}(C)$). There exists a line bundle $L$ on $C$ of degree
$2m+1$
such that:$$h^0\bigl( C,L(-c_1-\cdots -c_m)\bigr) >0\qquad{\rm and}\qquad
h^0\bigl(
C,L\otimes L^{-1}_a(-c_1-\cdots -c_m)\bigr) >0\ ,$$for all points
$c_1,\ldots,c_m$ of $C$.
Riemann-Roch implies that both $L$ and $L\otimes L^{-1}_a$ are $g^m_{2m+1}$'s.
By [ACGH],
exercise B-7, \page\ 138, this implies under our assumptions that $m=g-3$, and
that the
residual series of both $L$ and $L\otimes L^{-1}_a$ are $g^1_3$'s. If $g>4$,
there is at most one
$g^1_3$ and $a=0$. If $g=4$, there are at most two $g^1_3$'s and $a$ is either
$0$ or their
difference.\cqfd

\ind This proves that $\partial\CG {g+1,d}$ coincides over $\J g$ with
$\partial(\J {g+1}\cup\CT )$ for $g=4$ and $d=3$, and with
$\partial\J {g+1}$ otherwise. In particular:
{\nospacedmath $$\displaylines{\dim \partial{\cal F}\le 2+\dim \J
g=3g-1\qquad{\rm if}\quad
{\cal F}\supset \J {g+1}\ ,\cr
\dim \partial{\cal F}\le \dim \J 4=\dim\CT -1\qquad{\rm if}\quad {\cal
F}\supset
\CT\ .\cr}$$}\ind Since $\partial{\cal F}$ is the intersection of ${\cal F}$
with the
Cartier divisor $\partial\A {g+1}$, we have $
\dim {\cal F}\le\dim\partial{\cal F}+1$.
It follows that either ${\cal F}=\J {g+1}$ or ${\cal F}=\CT$, which finishes
the
proof of the induction step.

\ind To conclude, note that if $C$ is non-hyperelliptic of genus $4$, then
$g^1_3-h^1_3$ is
not in $C-C$ if it is non-zero. Therefore, $\CT$ is not contained in $\CG
{5,2}$.\cqfd

\ind Note that we only needed for $C-C$ to be a component of
$\partial\CG {g+1,d}\cap q^{-1}\{ (JC,\theta)\}$ for the proof. Our more
precise
result, combined with the fact that $\J 4=\CG {4,3}=\CG {4,2}$, proves:
$$\partial\CG {5,2}=\partial\J 5\qquad {\rm and}\qquad \partial\CG
{5,3}=\partial(\J
5\cup\CT)\ .$$
\ind This makes the following conjecture plausible, at least in dimension $5$:

{\pc Conjecture} -- For $1< d< g$ and $(g,d)\ne (5,3)$, then $\CG {g,d}=\J g$;
furthermore,\break $\CG {5,3}=\J 5\cup\CT$.

\ind A more tractable version would be the weaker:

{\pc Conjecture}' -- For $0< d< g$, then $\CG {g,d}\cap\CG {g,g-d}=\J g$.

\vfill\supereject\null
\saut \centerline {\pc References}

\saut
\hangindent=1cm
[ACGH] Arbarello, E., Cornalba, M., Griffiths, P., Harris, J.,
  {\it Geometry of algebraic curves. I. Grundlehren 267\/},
Springer-Verlag, New York, 1985.

\hangindent=1cm
[BC]	Barton, C., Clemens, C.H., {\it A result on the integral Chow ring of a
generic
principally polarized complex Abelian variety of dimension four\/}, Comp. Math.
{\bf 34} (1977), 49--67.

\hangindent=1cm
[B]	Beauville, A., {\it Sous-vari\'et\'es sp\'eciales des vari\'et\'es de
Prym\/}, Comp.
Math. {\bf 45} (1982), 357--383.

\hangindent=1cm
[CG]	Clemens, C.H., Griffiths, P.A., {\it The Intermediate Jacobian of the
Cubic
Threefold\/}, Ann. of Math. {\bf 95} (1972), 281--356.

\hangindent=1cm
[C] Collino, A., {\it A cheap proof of the irrationality of most cubic
threefolds\/}, Boll.
Un. Mat. Ital. {\bf 16} (1979), 451--465.

\hangindent=1cm
[HW]	Hulek, K., Weintraub, S., {\it The Principal Degenerations of Abelian
Surfaces and Their
Polarizations}, Math. Ann. {\bf 286} (1990), 281--307.

\hangindent=1cm
[M]	Matsusaka, T., {\it On a
characterization of a Jacobian Variety\/}, Mem. Coll. Sc. Kyoto, Ser. A, {\bf
23} (1959),
1--19.

\hangindent=1cm
[N]	Namikawa, Y., {\it A New Compactification of the Siegel Space and
Degeneration of Abelian
Varieties, I, II}, Math. Ann. {\bf 221} (1976), 97--141, 201--241.

\hangindent=1cm
[R1]	Ran, Z., {\it On subvarieties of abelian varieties\/}, Invent. Math. {\bf
62} (1981),
459--479.

\hangindent=1cm
[R2]	Ran, Z., {\it A Characterization of Five-Dimensional Jacobian
Varieties\/}, Invent. Math.
{\bf 67} (1982), 395--422.
\bigskip
\+ \kern 7cm &\cr
Olivier {\pc Debarre}
\+ URA D 0752 du CNRS & Department of Mathematics\cr
\+ Math\'ematique -- B\^atiment 425 & The University of Iowa\cr
\+ 91405 Orsay Cedex -- France & Iowa City  IA 52242 -- U.S.A.\cr
\+ debarre@matups.matups.fr & debarre@math.uiowa.edu \cr

\bye